\documentclass[12pt,a4paper,twoside]{article}

\newcommand{\pageformat}[6]{\setlength{\hoffset}{-1in}
                  \setlength{\voffset}{-1in}
                  \addtolength{\hoffset}{#5}
                            \addtolength{\voffset}{#6}
                            \setlength{\oddsidemargin}{#1}
                            \setlength{\evensidemargin}{#2}
                            \setlength{\textwidth}{\paperwidth}
                  \addtolength{\textwidth}{-\oddsidemargin}
                  \addtolength{\textwidth}{-\evensidemargin}
                  \addtolength{\textwidth}{-\marginparsep}
                  \addtolength{\textwidth}{-\marginparwidth}
                            \setlength{\topmargin}{#3}
                            \setlength{\textheight}{\paperheight}
                  \addtolength{\textheight}{-\topmargin}
                  \addtolength{\textheight}{-\headheight}
                  \addtolength{\textheight}{-\headsep}
                  \addtolength{\textheight}{-\footskip}
                  \addtolength{\textheight}{-#4}}
\pageformat{2cm}{3cm}{25mm}{25mm}{1pt}{0pt}

\usepackage{ifthen}
\newboolean{article}
    \setboolean{article}{true}
\newboolean{report}
\newboolean{book}
\newboolean{letter}
\newboolean{german}
\newboolean{italian}
\newboolean{nobaselinestretch}
\newboolean{nosectionappendix}
\newboolean{oldtoc}
\newboolean{nosectionequn}
\newboolean{notheorem}

\ifthenelse{\boolean{german}}{
    \usepackage{german}}{}

\usepackage[latin1]{inputenc}

\usepackage{amsmath}
\usepackage{amssymb}
\usepackage[mathscr]{eucal}

\ifthenelse{\boolean{notheorem}}{}{
    \usepackage{theorem}}

\ifthenelse{\boolean{nobaselinestretch}}{}{
    \renewcommand{\baselinestretch}{1.25}}

\newenvironment{env}[2]{\begin{#1}#2\end{#1}}{}
    \newcommand{\beq}[1]{\begin{env}{equation}{#1}}
    \newcommand{\beqn}[1]{\begin{env}{equation*}{#1}}
    \newcommand{\bal}[1]{\begin{env}{align}{#1}}
    \newcommand{\baln}[1]{\begin{env}{align*}{#1}}
    \newcommand{\bga}[1]{\begin{env}{gather}{#1}}
    \newcommand{\bgan}[1]{\begin{env}{gather*}{#1}}
    \newcommand{\bflal}[1]{\begin{env}{flalign}{#1}}
    \newcommand{\bflaln}[1]{\begin{env}{flalign*}{#1}}
    \newcommand{\bmu}[1]{\begin{env}{multline}{#1}}
    \newcommand{\bmun}[1]{\begin{env}{multline*}{#1}}
    \newcommand{\bsp}[1]{\begin{env}{split}{#1}}

    \newcommand{\eeq}{\end{env}}
    \newcommand{\eeqn}{\end{env}}
    \newcommand{\eal}{\end{env}}
    \newcommand{\ealn}{\end{env}}
    \newcommand{\ega}{\end{env}}
    \newcommand{\egan}{\end{env}}
    \newcommand{\eflal}{\end{env}}
    \newcommand{\eflaln}{\end{env}}
    \newcommand{\emu}{\end{env}}
    \newcommand{\emun}{\end{env}}
    \newcommand{\esp}{\end{env}}

\newcommand{\lf}{\vspace{2ex}}
\newcommand{\bulletline}[1][]{\lf\noindent~\hfill$\bullet\bullet\bullet$\hfill~

\lf\noindent\bf{#1}}

\renewcommand{\bf}[1]{\textbf{#1}}
\renewcommand{\it}[1]{\textit{#1}}

\renewcommand{\sf}[1]{\textsf{#1}}

\renewcommand{\tt}[1]{\texttt{#1}}
\newcommand{\hl}[1]{\bf{\it{#1}}}

\newcommand{\mbf}[1]{\mathbf{#1}}
\newcommand{\msf}[1]{\text{\small$\sf{#1}$}}

\newcommand{\cmc}[1]{\mathcal{#1}}
\newcommand{\eus}[1]{\mathscr{#1}}
\newcommand{\euf}[1]{\mathfrak{#1}}
\newcommand{\bb}[1]{\mathbb{#1}}

\newcommand{\nbd}[1]{$#1$\nobreakdash--}
\newcommand{\ol}[1]{\overline{#1}}

\newcommand{\wh}[1]{\widehat{#1}}
\newcommand{\ve}{\varepsilon}
\newcommand{\vt}{\vartheta}

\newcommand{\vp}{\varphi}
\newcommand{\om}{\omega}

\newcommand{\norm}[1]{\left\lVert#1\right\rVert}

\newcommand{\bnorm}[1]{\bigl\lVert#1\bigr\rVert}

\newcommand{\snorm}[1]{\norm{\smash{#1}}}

\newcommand{\bfam}[1]{\bigl(#1\bigr)}

\newcommand{\AB}[1]{\langle#1\rangle}

\newcommand{\CB}[1]{\{#1\}}
\newcommand{\bCB}[1]{\bigl\{#1\bigr\}}
\newcommand{\BCB}[1]{\Bigl\{#1\Bigr\}}
\newcommand{\SB}[1]{[#1]}

\newcommand{\RO}[1]{[#1)}

\newcommand{\set}[2][]{
    \ifthenelse{\equal{#1}{}}{
        \CB{#2}}{
        \CB{#1~|~#2}}}
\newcommand{\bset}[2][]{
    \ifthenelse{\equal{#1}{}}{
        \bCB{#2}}{
        \bCB{#1~|~#2}}}
\newcommand{\Bset}[2][]{
    \ifthenelse{\equal{#1}{}}{
        \BCB{#2}}{
        \BCB{#1~\big|~#2}}}

\DeclareMathOperator{\ls}{\normalfont\msf{span}}
\DeclareMathOperator{\cls}{\ol{\ls}}

\DeclareMathOperator{\id}{\normalfont\msf{id}}

\newcommand{\C}{\bb{C}}

\newcommand{\N}{\bb{N}}

\newcommand{\R}{\bb{R}}

\newcommand{\cA}{\cmc{A}}
\newcommand{\cB}{\cmc{B}}

\newcommand{\sB}{\eus{B}}

\newcommand{\sF}{\eus{F}}

\newcommand{\sK}{\eus{K}}

\newcommand{\sN}{\eus{N}}

\newcommand{\sS}{\eus{S}}

\newcommand{\eR}{\euf{R}}

\newcommand{\U}{\mbf{1}}

\ifthenelse{\boolean{nosectionequn}}{}{
    \numberwithin{equation}{section}
    }

\ifthenelse{\boolean{article}\or\boolean{letter}\or\boolean{nosectionequn}}{
    \setboolean{nosectionappendix}{true}}{}
\ifthenelse{\boolean{nosectionappendix}}{}{
    \renewcommand{\appendix}{
        \chapter*{\appendixname}
        \addcontentsline{toc}{chapter}{\appendixname}
        \renewcommand{\thesection}{\Alph{section}}
        \setcounter{section}{0}}}
   
\ifthenelse{\boolean{report}\or\boolean{book}}{
    }{}

\ifthenelse{\boolean{notheorem}}{}{
        \newcommand{\mnname}{Mathematical note.}
        \newcommand{\enname}{End of the note.}
        \newcommand{\definame}{Definition.}
        \newcommand{\propname}{Proposition.}
        \newcommand{\lemname}{Lemma.}
        \newcommand{\exname}{Example.}
        \newcommand{\exername}{Exercise.}
        \newcommand{\remname}{Remark.}
        \newcommand{\obname}{Observation.}
        \newcommand{\thmname}{Theorem.}
        \newcommand{\corname}{Corollary.}
        \newcommand{\proofname}{Proof.}
    \ifthenelse{\boolean{german}}{
        \renewcommand{\mnname}{Mathematische Notiz.}
        \renewcommand{\enname}{Ende der Notiz.}
        \renewcommand{\exname}{Beispiel.}
        \renewcommand{\exername}{Übung.}
        \renewcommand{\remname}{Bemerkung.}
        \renewcommand{\obname}{Beobachtung.}
        \renewcommand{\thmname}{Satz.}
        \renewcommand{\corname}{Korollar.}
        \renewcommand{\proofname}{Beweis.}}{}
    \ifthenelse{\boolean{italian}}{
        \renewcommand{\mnname}{Nota matematica.}
        \renewcommand{\enname}{Fina della nota.}
        \renewcommand{\definame}{Definizione.}
        \renewcommand{\propname}{Proposizione.}
        \renewcommand{\exname}{Esempio.}
        \renewcommand{\exername}{Esercizio.}
        \renewcommand{\remname}{Nota.}
        \renewcommand{\obname}{Osservazione.}
        \renewcommand{\thmname}{Teorema.}
        \renewcommand{\corname}{Corollario.}
        \renewcommand{\proofname}{Dimostrazione.}

       \renewcommand{\appendixname}{Appendice}

       }{}
    \theoremheaderfont{\normalfont\bfseries}
    \theoremstyle{change}
        \theorembodyfont{\rmfamily}
            \newtheorem{emp}{}[section]
                \newcommand{\bemp}[1][]{
                    \begin{emp}\hskip-\labelsep\bf{#1}\hskip\labelsep}
                \newcommand{\eemp}{\end{emp}}
\newtheorem{itemp}[emp]{}
                \newcommand{\bitemp}[1][]{
                    \begin{itemp}\hskip-\labelsep\bf{#1}\hskip\labelsep\normalfont\itshape}
                \newcommand{\eitemp}{\end{itemp}}
            \newtheorem{mn}[emp]{\mnname}
                \newcommand{\bnm}{\begin{mn}~\begin{quotation}\renewcommand{\baselinestretch}{1}\small\noindent\ignorespaces}
                \newcommand{\enm}{\end{quotation}\hfill\bf{\enname}\end{mn}}
            \newtheorem{ex}[emp]{\exname}
                \newcommand{\bex}{\begin{ex}}
                \newcommand{\eex}{\end{ex}}
            \newtheorem{exer}[emp]{\exername}
                \newcommand{\bexer}{\begin{exer}}
                \newcommand{\eexer}{\end{exer}}
            \newtheorem{defi}[emp]{\definame}
                \newcommand{\bdefi}{\begin{defi}}
                \newcommand{\edefi}{\end{defi}}
            \newtheorem{rem}[emp]{\remname}
                \newcommand{\brem}{\begin{rem}}
                \newcommand{\erem}{\end{rem}}
            \newtheorem{ob}[emp]{\obname}
                \newcommand{\bob}{\begin{ob}}
                \newcommand{\eob}{\end{ob}}
        \theorembodyfont{\normalfont\itshape}
            \newtheorem{thm}[emp]{\thmname}
                \newcommand{\bthm}{\begin{thm}}
                \newcommand{\ethm}{\end{thm}}
            \newtheorem{prop}[emp]{\propname}
                \newcommand{\bprop}{\begin{prop}}
                \newcommand{\eprop}{\end{prop}}
            \newtheorem{cor}[emp]{\corname}
                \newcommand{\bcor}{\begin{cor}}
                \newcommand{\ecor}{\end{cor}}
            \newtheorem{lem}[emp]{\lemname}
                \newcommand{\blem}{\begin{lem}}
                \newcommand{\elem}{\end{lem}}
\newenvironment{empn}[1]{\lf\noindent\bf{#1}\ignorespaces\hskip\labelsep}{\lf}
		\newcommand{\bempn}[1]{\begin{empn}{#1}}
		\newcommand{\eempn}{\end{empn}}
		\newcommand{\bitempn}[1]{\begin{empn}{#1}\normalfont\itshape}
		\newcommand{\eitempn}{\end{empn}}
                \newcommand{\bnmn}{\begin{empn}{\mnname}~\begin{quotation}\renewcommand{\baselinestretch}{1}\small\noindent\ignorespaces}
                \newcommand{\enmn}{\end{quotation}\hfill\bf{\enname}\end{empn}}
		\newcommand{\bexn}{\begin{empn}{\exname}}
		\newcommand{\eexn}{\end{empn}}
		\newcommand{\bexern}{\begin{empn}{\exername}}
		\newcommand{\eexern}{\end{empn}}
		\newcommand{\bdefin}{\begin{empn}{\definame}}
		\newcommand{\edefin}{\end{empn}}
		\newcommand{\bremn}{\begin{empn}{\remname}}
		\newcommand{\eremn}{\end{empn}}
		\newcommand{\bobn}{\begin{empn}{\obname}}
		\newcommand{\eobn}{\end{empn}}

\newcommand{\qedsymbol}{~\rule[-0.35mm]{2mm}{2mm}}
    \newcounter{proof}[emp]
    \newenvironment{Proof}[1]{
        \vspace{1ex}
        \renewcommand{\item}[1][\stepcounter{proof}(\roman{proof})]%
            {##1\hskip\labelsep}
        \noindent\textsc{#1\hskip\labelsep}}{
        \nolinebreak\qedsymbol}
    \newcommand{\proof}[1][\proofname]{
        \begin{Proof}{#1}\ignorespaces}
    \newcommand{\qed}{\end{Proof}}
    \newcommand{\noqed}{
        \renewcommand{\qedsymbol}{}
        \end{Proof}}}
    \ifthenelse{\boolean{italian}}{
        \renewcommand{\proofname}{Dimostrazione.}}{}

\usepackage[varg]{txfonts}

\usepackage[hypertex]{hyperref}

\begin{document}

\bibliographystyle{amsalpha}

\title{Nondegenerate Representations of \\Continuous Product Systems\thanks{AMS 2000 subject classification 46L55, 46L08, 46L53, 60G20}}
\author{}
\author{
~\\
Michael Skeide\thanks{This work is supported by research funds of University of Molise and Italian MIUR.}\\\\
{\small\itshape Universit\`a\ degli Studi del Molise}\\
{\small\itshape Dipartimento S.E.G.e S.}\\
{\small\itshape Via de Sanctis}\\
{\small\itshape 86100 Campobasso, Italy}\\
{\small{\itshape E-mail: \tt{skeide@unimol.it}}}\\
{\small{\itshape Homepage: \tt{http://www.math.tu-cottbus.de/INSTITUT/lswas/\_skeide.html}}}\\
}
\date{}

{
\renewcommand{\baselinestretch}{1}
\maketitle



\begin{abstract}
\noindent
We show that every (continuous) faithful product system admits a (continuous) faithful nondegenerate representation. For Hilbert spaces this is equivalent to Arveson's result that every Arveson system comes from an \nbd{E_0}semigroup. We point out that for Hilbert modules this is not so. As applications we show a \nbd{C^*}algebra version of a result for von Neumann algebras due to Arveson and Kishimoto, and a result about existence of elementary dilations for (semi-)faithful CP-semigroups.
\end{abstract}

}



\section{Introduction}

Recall that a \hl{correspondence} over $\cB$ (or, more generally, from $\cA$ to $\cB$) is a (right) Hilbert \nbd{\cB}module with a \hl{nondegenerate} left action of $\cB$ (or of $\cA$) by adjointable operators. The (internal) tensor product of correspondences we denote by $\odot$. The \hl{algebra of adjointable operators} on a Hilbert \nbd{\cB}module we denote by $\sB^a(E)$. Recall that $E$ is \hl{full}, if $\cls\AB{E,E}=\cB$. We say a correspondence is \hl{faithful}, if its left action defines a faithful representation. For $x\in E$, we define the mapping $x^*\colon E\rightarrow\cB$ by setting $x^*y=\AB{x,y}$. The adjoint of $x^*$ is $x\colon b\mapsto xb$. The algebra $\sK(E)$ of \hl{compact operators} on $E$ is the norm completion of the algebra $\sF(E):=\ls\CB{xy^*\colon x,y\in E}$ of \hl{finite rank operators} that is spanned linearly by the \hl{rank-one operators} $xy^*\colon z\mapsto x\AB{y,z}$. Recall that for a unital homomorphism $\vt\colon\sB^a(E)\rightarrow\sB^a(F)$ to be \hl{strict} is equivalent to that the action of the compacts alone is already nondegenerate: $\cls\vt(\sK(E))F=F$.

An (algebraic) \hl{product system} is a family $E^\odot=\bfam{E_t}_{t\in\R_+}$ of correspondences $E_t$ over a \nbd{C^*}al\-gebra $\cB$ with a family of bilinear unitaries $u_{s,t}\colon E_s\odot E_t\rightarrow E_{s+t}$, such that the ``multiplication'' defined by $x_sy_t:=u_{s,t}(x_s\odot y_t)$ is associative. Moreover, $E_0=\cB$ is the trivial correspondence over $\cB$ and $u_{0,t}$ and $u_{t,0}$ are left and right action, respectively, of $E_0=\cB$ on $E_t$. A product system $E^\odot$ is \hl{full}, if each $E_t$ $(t\ge0)$ is full, and $E^\odot$ is \hl{faithful}, if each $E_t$  $(t\ge0)$ is faithful.

A \hl{left dilation} of a full product system $E^\odot$ to a full Hilbert \nbd{\cB}module $E$ is a family of unitaries $v_t\colon E\odot E_t\rightarrow E$ such that $(xy_s)z_t=x(y_sz_t)$, where we defined $xy_t:=v_t(x\odot y_t)$. By setting $\vt^v_t(a):=v_t(a\odot\id_t)v_t^*$, every left dilation gives rise to an \hl{\nbd{E_0}semigroup} $\vt^v=\bfam{\vt^v_t}_{t\in\R_+}$ on $\sB^a(E)$, that is, a semigroup of strict unital endomorphisms. (In these notes, as a convention, we always assume that homomorphisms are strict.)

A \hl{right dilation} of a faithful product system $E^\odot$ to a faithful correspondence $H$ from $\cB$ to $\C$ (that is, a Hilbert space with a faithful nondegenerate representation of $\cB$) is a family of bilinear unitaries $w_t\colon E_t\odot H\rightarrow H$ such that $x_t(y_sh)=(x_ty_s)h$, where we defined $x_th:=w_t(x_t\odot h)$. By defining the operator $\eta^w_t(x_t)\colon h\mapsto x_th$ in $\sB(H)$, every right dilation gives rise to a \hl{representation} of $E^\odot$, that is, a family of linear maps $\eta^w_t\colon E_t\rightarrow\sB(H)$ such that
\beqn{
\eta^w_t(x_t)^*\eta^w_t(y_t)
~=~
\eta^w_0(\AB{x_t,y_t})
\text{~~~~~~and~~~~~~}
\eta^w_t(x_t)\eta^w_s(y_s)
~=~
\eta^w_{t+s}(x_ty_s),
}\eeqn
which is \hl{nondegenerate} (that is, $\cls\eta^w_t(E_t)H=H$ for all $t$) and \hl{faithful} (that is, $\eta^w_0$ and, therefore, all $\eta^w_t$ are injective). By applying the two equations to $\AB{x_0,y_0}=x_0^*y_0$ $(x_0,x_0,^*,y_0\in E_0=\cB)$, and taking also into account nondegeneracy, one checks that $\eta^w_0$ is a representation of $\cB$. (Note, too, that $\eta^w_t(x_t\AB{y_t,z_t})=\eta^w_t(x_t)\eta^w_t(y_t)^*\eta^w_t(z_t)$, that is, $\eta_t$ is a \hl{ternary homomorphism}. In particular, $\eta^w_t$ is linear and completely contractive; see Abbaspour and Skeide \cite{AbSk07}.) Conversely, if $\eta_t$ is a faithful nondegenerate representation of $E^\odot$ on $H$, then $H$ is a faithful correspondence from $\cB$ to $\C$ via $\eta_0$ and by setting $w_t(x_t\odot h):=\eta_t(x_t)h$ we define a right dilation to $H$. Of course, $w_t$ gives back $\eta_t$ as $\eta^w_t$.

\lf
We see that left dilations relate full product systems to \nbd{E_0}semigroups, while right dilations of faithful product systems are synonymous with nondegenerate faithful representations.

An \nbd{E_0}semigroup $\vt$ on $\sB^a(E)$ with $E$ a full Hilbert \nbd{\cB}module give rise to a full product system $E^\odot$ of \nbd{\cB}correspondences and a left dilation $v_t$ such that $\vt=\vt^v$; see \cite{Ske02,Ske09,Ske09a}. Two \nbd{E_0}semigroups on the same $\sB^a(E)$ have isomorphic product systems if an only if they are cocycle conjugate; see \cite{Ske02}. In \cite{Ske08p1,Ske09a} we extended this to \nbd{E_0}semigroups acting on different $\sB^a(E)$ provided the two $E$ are countably generated and over unital $\cB$ (\cite{Ske08p1}) or \nbd{\sigma}unital $\cB$ (\cite{Ske09a}). In Skeide \cite{Ske07}, we have constructed for every \it{continuous} product system of correspondences over a \it{unital} \nbd{C^*}algebra $\cB$ a \it{continuous} left dilation, that is, it is the product system of a strongly continuous \nbd{E_0}semigroup. This dilation is to a countably generated $E$, if the product system is \it{countably generated}. Also this we generalized to \nbd{\sigma}unital $\cB$ in \cite{Ske09a}.\footnote{We use the occasion to mention that the proof in \cite[Proposition 4.9]{Ske07} that the product system of the constructed \nbd{E_0}semigroup has the same continuous structure as the original product system, has a gap. This gap is fixed in \cite{Ske09a}.} Combining all this, under the stated countability assumptions in \cite{Ske08p1} we obtain a full analogy with Arveson's results \cite{Arv89,Arv90a,Arv89a,Arv90}, namely, a one-to-one correspondence between \nbd{E_0}semigroups on full Hilbert \nbd{\cB}modules up to stable cocycle conjugacy and full product systems up to isomorphism.

In these notes we show that every (continuous) faithful product system of correspondences over an arbitrary \nbd{C^*}algebra admits a (continuous) right dilation, that is, a (continuous) faithful nondegenerate representation (Theorem \ref{rdilthm}). Right dilations do not establish such a direct relation between product systems and \nbd{E_0}semigroups. It is true that a right dilation $\bfam{w_t}_{t\in\R_+}$ of $E^\odot$ gives rise to an \nbd{E_0}semigroup, namely, to the semigroup of unital normal endomorphisms $\theta^w_t(a):=w_t(\id_t\odot a)w_t^*$ acting on $\sB^{bil}(H)$, the von Neumann subalgebra of $\sB(H)$ which consists of all bounded bilinear operators on $H$. However, the \nbd{E_0}semigroup $\theta^w$, in general, does not allow to reconstruct the product system $E^\odot$ uniquely. The following space of intertwining operators $\bCB{x_t\in\sB(H)\colon\theta_t^w(a)x_t=x_ta~(a\in\sB^{bil}(H))}$ contains $E_t$ as a (strongly dense) subset in a natural way. But $E_t$ will coincide with that intertwiner space, only if it is a von Neumann correspondence over the double commutant of $\eta^w_0(\cB)$ in $\sB(H)$. (We omit details and refer the reader to \cite{Ske03c,Ske08,Ske06b}.) Nevertheless, the question whether a product system admits a faithful nondegenerate representation is of independent interest. In Section \ref{applSEC} we give some applications. We prove a result about \it{embedding} faithful \nbd{E_0}semigroups into inner automorphism groups (which provides an analogue for \nbd{C^*}algebras of a result by Arveson and Kishimoto for von Neumann algebras). And we prove existence of \it{elementary dilations} for (semi-)faithful CP-semigroups.

Technically, these notes where we show that every faithful product system admits a right dilation, is very similar to Skeide \cite{Ske07}, where we constructed a left dilation for every \it{continuous} product system of correspondences over a \it{unital} \nbd{C^*}algebra. In cases where the results are just analogues of statements in \cite{Ske07} with analogue proofs (Propositions \ref{rcllprop} -- \ref{cdprop}), we do not repeat these proofs. Proposition \ref{Hcontprop} and its corollary, instead, are technically more involved (mainly, because we consider sections in spaces that are tensor products), and require new ideas.

In principle, the reader who is interested only in the proof of the statement, may now pass immediately to Theorem \ref{zeta'thm} and, then, proceed to Section \ref{constrSEC}. But we wish to clarify, why such a strange statement like Theorem \ref{zeta'thm} is the natural starting point for the construction of a right dilation, knowing the successful strategy for constructing a left dilation. The balance of this introduction is dedicated to this motivation.

\bulletline
For product systems of Hilbert spaces $E^\otimes$ (\hl{Arveson systems}) there is not much a difference between the construction of a left dilation and the construction of a right dilation. More precisely, a left dilation of $E^\otimes$ gives rise to a right dilation of ${E^\otimes}^{op}$ (the \hl{opposite} Arveson system of $E^\otimes$ with the opposite product $(x_s,y_t)\mapsto y_tx_s$) and \it{vice versa}, simply by ``inverting'' all orders in tensor products; see \cite{Ske06}. (Note, however, that $E^\otimes$ and ${E^\otimes}^{op}$ need not be isomorphic Arveson systems; see Tsirelson \cite{Tsi00p1}.)

For Hilbert modules the situation is more delicate. Since there is no canonical flip operation for tensor products of correspondences, there is no such thing like the opposite product system of $E^\odot$. However, there is the \it{commutant} of von Neumann correspondences (Skeide \cite{Ske03c} and Muhly and Solel \cite{MuSo04}). The commutant transforms a product system of von Neumann correspondences over a von Neumann algebra $\cB\subset\sB(G)$ into a product system of von Neumann correspondences over the commutant $\cB'$ of $\cB$. (In fact, the opposite of an Arveson system is just its commutant system.) Under commutant, left dilations transform into right dilations and \it{vice versa}; \cite[Theorem 9.9]{Ske09} or \cite[Theorem 3.6(3)]{Ske08}. Also the conditions to be full or faithful are interchanged under commutant.

We see, in the von Neumann case, a proof of existence of left dilation transforms into a proof of existence of right dilations. Although, this is strictly true only for von Neumann correspondences, understanding how ingredients of the proof for existence of left dilations transform under the commutant, is crucial for finding our proof here for existence of right dilations in the case of \nbd{C^*}correspondences.

For instance, in Skeide \cite{Ske06} we constructed a left and right dilation for every Arveson $E^\otimes=\bfam{E_t}_{t\in\R_+}$ system by starting with a left and right dilation, respectively, of the discrete subsystem $\bfam{E_t}_{t\in\N_0}$ of $E^\otimes$, and ``blowing it up'' suitably. As already mentioned, there is not really a difference between left and right, here. It is, however, important to note that the input, a dilation of the discrete subsystem, can easily be obtained by choosing a unit vector $\zeta_1\in E_1$.

In \cite{Arv06}, Arveson constructed a right dilation, which turned out to be unitarily equivalent to ours; see Skeide \cite{Ske06a}. In order to construct that right dilation, also Arveson fixes a unit vector $\zeta_1\in E_1$. Then he considers the space of \it{right stable} sections in $L^{2,loc}(E^\otimes)$, that is, of locally square integrable sections $x\colon\alpha\mapsto x_\alpha\in E_\alpha$, that fulfill $x_{\alpha+1}=x_\alpha\zeta_1$ for all sufficiently big $\alpha$. He equips this space with a semiinner product $\AB{x,y}:=\int_T^{T+1}\AB{x_\alpha,y_\alpha}\,d\alpha$ (which does not depend on $T$ for all sufficiently big $T$). On the quotient Hilbert space the product system acts simply by ``multiplication'' (that is, tensor product) from the left.

The same construction works for left dilations. We simply start with the space of \it{left stable} sections (that is, $x_{\alpha+1}=\zeta_1x_\alpha$ for all sufficiently big $\alpha$), on which the product system acts by ``multiplication'' from the right. This construction of a left dilation also works for Hilbert modules as soon as we have a \hl{unit vector} $\zeta_1\in E_1$ (that is, a vector that fulfills $\AB{\zeta_1,\zeta_1}=\U\in\cB$), because then the semiinner product $\int_T^{T+1}\AB{x_\alpha,y_\alpha}\,d\alpha$ does not depend on $T$ for all sufficiently big $T$. Note that $\cB$ must be unital. \it{Continuous} product systems of correspondences over a unital \nbd{C^*}algebras have unit vectors in all fibres, and apart from more involved technical problems due to modules, the proof in \cite{Ske07} runs like Arveson's \cite{Arv06}.

Nothing like this is possible for right dilations in the module case: Since $\AB{x_\alpha\odot\zeta_1,y_\alpha\odot\zeta_1}=\AB{\zeta_1,\AB{x_\alpha,y_\alpha}\zeta_1}$, the semiinner product $\int_T^{T+1}\AB{x_\alpha,y_\alpha}\,d\alpha$ of right stable sections will depend on $T$, unless $\zeta_1$ is a unit vector that commutes with all elements in $\cB$.

The solution to our becomes clear, if we recall that the unit vector corresponds to the construction of a certain left dilation of the discrete subsystem $\bfam{E_t}_{t\in\N_0}$ of $E^\odot$. In the construction of a right dilation, dual to the construction of a left dilation, that unit vector must be a unit vector in the commutant $E'_1$ of $E_1$, rather than in $E_1$.

The construction of the commutant of a correspondence $E_1$ over $\cB$ is possible, as soon as we assume that $\cB\subset\sB(G)$ is a concrete \nbd{C^*}algebra of operators acting nondegenerately on a Hilbert space $G$: Define the correspondence $H_1:=E_1\odot G$ from $\cB$ to $\C$. Put $E'_1:=\sB^{bil}(G,H_1)$.

\brem
On $H_1$ we have an action of the von Neumann algebra $\cB'=\sB^{bil}(G)$ defined by $b'(x_1\odot g):=x_1\odot b'g$, the so-called \hl{commutant lifting}.  $E'_1$ with inner product $\AB{x'_1,y'_1}:=x'^*_1y'_1$ and with left action via the commutant lifting is, then, a von Neumann \nbd{\cB'}correspondence. If $\cB\subset\sB(G)$ is a von Neumann algebra and if $E_1$ is a von Neumann \nbd{\cB}correspondence, then $E'_1$ is precisely the \hl{commutant} of $E_1$ as introduced in \cite{Ske03c}.
\erem

We see that a unit vector in $E'_1$ is an isometry $\zeta'_1$ from $G$ to $H_1$ that intertwines the canonical actions of $\cB$.

Existence of an identification $\cB\subset\sB(G)$ such that there exists a unit vector $\zeta'_1\in E'_1$ follows now by existence of a nondegenerate faithful representation of $E_1$ as proved by \cite{Hir05a,Ske09}.

\bthm\label{zeta'thm}
Let $\cB$ be a \nbd{C^*}algebra and $E$ a faithful correspondence over $\cB$. Then there exists a faithful nondegenerate representation of $\cB$ on a Hilbert space $G$ that admits an isometry $\zeta'\in\sB^{bil}(G,E\odot G)$.
\ethm

\proof
Under the same hypothesis, in \cite[Theorem 8.3]{Ske09} (first \cite{Hir05a} for the case when $E$ is also full) we have shown that $E$ admits a \hl{faithful nondegenerate representation}. That is, there exists a faithful nondegenerate representation $\pi\colon\cB\rightarrow\sB(G)$ and a map $\eta\colon E\rightarrow\sB(G)$ such that
\baln{
\eta(b_1xb_2)
&
~=~
\pi(b_1)\eta(x)\pi(b_2),
&
\pi(\AB{x,y})
&
~=~
\eta(x)^*\eta(y),
&
\cls\eta(E)G
&
~=~
G.
}\ealn
We immediately check that $E\odot G=G$ as \nbd{\cB}\nbd{\C}correspondences, via $x\odot g\mapsto\eta(x)g$. Clearly, $\zeta':=\id_G$ is an isometry in $\sB^{bil}(G)=\cB'$.\qed

\section{The construction}\label{constrSEC}

In principle, like all the constructions in \cite{Ske06,Arv06,Ske07}, everything works al\-ge\-braical\-ly, if we take the occurring direct integrals with respect to the counting measure, that is, if we take direct sums. But nice continuity (or just measurability) properties of the constructed dilations, of course, emerge only if we choose the Lebesgue measure. For that goal the product system must fulfill technical conditions. Here, as in \cite{Ske07}, we consider \it{continuous} product systems as defined in \cite{Ske03b}. Also measurable versions of product systems have been considered; see Hirshberg \cite{Hir04}.

\bdefi\label{cPSdef}
Let $E^\odot=\bfam{E_t}_{t\in\R_+}$ be a product system of correspondences over a \nbd{C^*}algebra $\cB$ with a family $i=\bfam{i_t}_{t\in\R_+}$ of \hl{isometric} (that is, inner product preserving) embeddings $i_t\colon E_t\rightarrow\wh{E}$ into a Hilbert \nbd{\cB}module $\wh{E}$. Denote by
$$
CS_i(E^\odot)
~=~
\BCB{x=\bfam{x_t}_{t\in\R_+}\colon x_t\in E_t,t\mapsto i_tx_t\mbox{~is continuous}}
$$
the set of \hl{continuous sections} of $E^\odot$ (with respect to $i$). We say $E^\odot$ is \hl{continuous} (with respect to $i$), if the following conditions are satisfied.
\begin{enumerate}
\item\label{csdef}
For every $y_t\in E_t$ we can find a continuous section $x\in CS_i(E^\odot)$ such that $x_t=y_t$.

\item \label{cpsdef}
For every pair $x,y\in CS_i(E^\odot)$ of continuous sections the function
$$
(s,t)
~\longmapsto~
i_{s+t}(x_sy_t)
$$
is continuous.
\end{enumerate}
We say two embeddings $i$ and $i'$ have the same \hl{continuous structure}, if $CS_i(E^\odot)=CS_{i'}(E^\odot)$.
\edefi

This is \cite[Definition 7.1]{Ske03b} except that $\cB$ need not be unital. It is motivated by the fact that every product system of a strictly continuous \nbd{E_0}semigroup acting on the operators of a Hilbert module fulfills these requirements. Condition 1 may be replaced with the weaker condition that for every $t\in\R_+$ the set $\CB{x_t\colon x\in CS_i(E^\odot)}$ is total in $E_t$. Condition 2 may be replaced with the weaker condition that for all $b\in\cB$ with $x\in CS_i(E^\odot)$ also the section $bx:=\bfam{bx_t}_{t\in\R_+}$ is in $CS_i(E^\odot)$ and that the function $(s,t)\mapsto \AB{z,i_{s+t}(x_sy_t)}$ is continuous for every $z\in\wh{E}$ and every pair $x,y\in CS_i(E^\odot)$. See \cite{Ske03b,Ske07} for details.

\lf
Throughout the balance of this section we shall suppose that $E^\odot$ is faithful and that $\cB\subset\sB(G)$ is given as a concrete \nbd{C^*}algebra of operators via the representation guaranteed by Theorem \ref{zeta'thm} for the correspondence $E=E_1$.

For each $t\in\R_+$ we put $H_t:=E_t\odot G$. For all $s,t\in\R_+$ an element $x_s\in E_s$ has an action $h_t\mapsto x_sh_t\in H_{s+t}$ on $h_t\in H_t$ defined by setting $x_s(y_t\odot g):=(x_sy_t)\odot g$.

Each $H_t$ is a correspondence from $\cB$ to $\C$ (faithful, if and only if $E_t$ is faithful) with left action defined by setting $b(x_t\odot g):=(bx_t)\odot g$.

By Theorem \ref{zeta'thm}, we may fix an isometry $\zeta'_1\in\sB^{bil}(G,H_1)$. It would be tempting to  consider \it{stable} sections $x$ of $E^\odot$, in the sense that $x_{\alpha+1}=x_\alpha\zeta'_1$ for all sufficiently large $\alpha$. However, if $x_\alpha$ is in $E_\alpha$, then $x_\alpha\zeta'_1\in\sB(G,H_{\alpha+1})$ is, in general, not in $E_{\alpha+1}=\bset[g\mapsto x_{\alpha+1}\odot g]{x_{\alpha+1}\in E_{\alpha+1}}\subset\sB(G,H_{\alpha+1})$. Thus, we cannot consider sections in $E^\odot$. Instead, we will consider sections in the family $E^\odot\odot G:=\bfam{H_t}_{t\in\R_+}$ of Hilbert spaces. As $\zeta'_1$ is left linear, for each $t$ the operator $\id_t\odot\zeta'_1$ is a well-defined isometry in $\sB^{bil}(H_t,E_t\odot H_1)$. As $E_t\odot H_1\cong H_{t+1}$ (as correspondences from $\cB$ to $\C$) via $x_t\odot h_1\mapsto x_th_1$, there is also a unique isometry in $\sB^{bil}(H_t,H_{t+1})$ that sends $x_t\odot g$ to $x_t\zeta'_1g$. Also here we shall write the action of this isometry on elements $h_t\in H_t$ simply as $h_t\mapsto\zeta'_1h_t$.

\brem
This notation suggests a relation $\zeta'_1x_s=x_s\zeta'_1$ for the actions on $H_t$. A fact, that may be verified for all pairs $x_s\in E_s$ and $y'_r\in E'_r:=\sB^{bil}(G,H_r)$; see the proof of \cite[Theorem 3.4(1)]{Ske08}.
\erem

The sections $h$ of $E^\odot\odot G$ we shall consider, will fulfill
\beqn{\tag{$**$}\label{stabsec}
h_{\alpha+1}
~=~
\zeta'_1h_\alpha
}\eeqn
for all sufficiently big $\alpha$. But, before we proceed we need to specify some properties of the relevant direct integrals.

First of all, we note that the embeddings $i_t\colon E_t\rightarrow\wh{E}$ give rise to embeddings $H_t\rightarrow\wh{H}:=\wh{E}\odot G$, also denoted by $i_t$, defined by $x_t\odot g\mapsto (i_tx_t)\odot g$. We, therefore, may speak about continuous sections $h=\bfam{h_t}_{t\in\R_+}$ of $E^\odot\odot G$, in the sense that $t\mapsto i_th_t$ is continuous. We denote the set of all continuous sections of $E^\odot\odot G$ by $CS_i(E^\odot\odot G)$. Obviously, whenever $x,y\in CS_i(E^\odot)$ and $g\in G$, then the functions $t\mapsto i_tx_t\odot g$ and $(s,t)\mapsto i_{s+t}(x_sy_t)\odot g$ are continuous. The following properties are less obvious.

\bprop\label{Hcontprop}
Every continuous section $h\in CS_i(E^\odot\odot G)$ may be approximated locally uniformly by elements in $\ls CS_i(E^\odot)\odot G$. Moreover:
\begin{enumerate}
\item\label{Hcsdef}
For every $k_t\in H_t$ we can find a continuous section $h\in CS_i(E^\odot\odot G)$ such that $h_t=k_t$.

\item \label{Hcpsdef}
For every pair $x\in CS_i(E^\odot)$ and $h\in CS_i(E^\odot\odot G)$ of continuous sections the function
$$
(s,t)
~\longmapsto~
i_{s+t}(x_sh_t)
$$
is continuous.
\end{enumerate}
\eprop

\proof
The proof of \eqref{Hcsdef} is very similar to the proof of \cite[Proposition 7.9]{Ske03b}. Every $k_t\in H_t$ may be written as $\sum_nk^n_t$ in such a way that $\sum_n\snorm{k^n_t}<\infty$, where $k^n_t=\sum_{j=1}^{m_n}y^{n,j}_t\odot f^{n,j}$ with $y^{n,j}_t\in E_t$, $f^{n,j}\in G$ $(n\in\N;m_n\in\N;j=1,\ldots,m_n)$. Choose continuous sections $x^{n,j}\in CS_i(E^\odot)$ such that $x^{n,j}_t=y^{n,j}_t$ and such that $\bnorm{\sum_{j=1}^{m_n}x^{n,j}_s\odot f^{n,j}}\le\snorm{k^n_t}$ for all $s\in\R_+$, $n\in\N$. Then $h:=\sum_n\sum_{j=1}^{m_n}x^{n,j}\odot f^{n,j}$ is a (bounded!) section in $CS_i(E^\odot\odot G)$ with $h_t=k_t$.

Now let $h\in CS_i(E^\odot\odot G)$ and choose $0\le a<b<\infty$ and $\ve>0$. For every $\beta\in\SB{a,b}$, by the proof of Part \eqref{Hcsdef} there exists a section $h^\beta$ in $\ls CS_i(E^\odot)\odot G$ such that $\snorm{h_\beta-h^\beta_\beta}<\frac{\ve}{2}$. For every $\beta$ define $I_\beta$ to be the largest interval such that $\snorm{h_\alpha-h^\beta_\alpha}<\ve$ for all $\alpha\in I_\beta$. Every $I_\beta$ is open in $\SB{a,b}$ and contains at least $\beta$. Therefore, the family of all $I_\beta$ forms an open cover of the compact interval $\SB{a,b}$. So, we may choose $\beta_1,\ldots,\beta_m$ such that the union over $I_{\beta_i}$ is $\SB{a,b}$. By standard theorems about \it{partitions of unity} there exist continuous functions $\vp_i$ on $\SB{a,b}$ with the following properties:
\baln{
0
~\le~
&
\vp_i
~\le~
1,
&
\vp_i\upharpoonright I_{\beta_i}^\complement
&
~=~
0,
&
\sum_{i=1}^m\vp_i
&
~=~
1.
}\ealn
From these properties, one easily verifies that $\norm{h_\alpha-\sum_{i=1}^m\vp_i(\alpha)h^{\beta_i}_\alpha}<\ve$ for all $\alpha\in\SB{a,b}$. This shows that $\sum_{i=1}^m\vp_ih^{\beta_i}\in\ls CS_i(E^\odot)\odot G$ approximates $h$ uniformly up to $\ve$ on the interval $\SB{a,b}$.

\eqref{Hcpsdef} follows now by three epsilons, approximating $h_t$ with an element in $\ls CS_i(E^\odot)\odot G$ on a suitably big interval.\qed

\bcor\label{shiftcor}
If $h\colon t\mapsto h_t$ is a continuous section, then the shifted section
\beqn{
t
~\longmapsto~
\begin{cases}
0&t<1
\\
\zeta'_1h_{t-1}&t\ge1
\end{cases}
}\eeqn
is continuous for $t\ge1$.
\ecor

\proof
By Proposition \ref{Hcontprop} the elements in $\ls CS_i(E^\odot)\odot G$ approximate $t\mapsto h_t$ locally uniformly. So, it is enough to show the statement for sections of the form $t\mapsto x_t\odot g$ $(x\in CS_i(E^\odot),g\in G)$. Again by Proposition \ref{Hcontprop} there is a section $h\in CS_i(E^\odot\odot G)$ such that $h_1=\zeta'_1g$. Once more, by Proposition \ref{Hcontprop} we may choose $h_t=\sum_{n=1}^\infty y_t^n\odot g_n$ (limit locally uniformly in $t$) for sequences $y^n\in CS_i(E^\odot)$ and $g^n\in G$. Then
\beqn{
\zeta'_1(x_t\odot g)
~=~
x_t\zeta'_1g
~=~
\sum_{n=1}^\infty x_ty_1^n\odot g_n
}\eeqn
locally uniformly in $t\ge0$. All $t\mapsto x_ty_1^n\odot g_n$ are continuous locally uniformly, so that also $t\mapsto \zeta'_1(x_t\odot g)$ is locally uniformly continuous.\qed

\lf
In the sequel, for a section $h$ of $E^\odot\odot G$ we shall denote $h(t):=i_th_t$. Let $0\le a<b<\infty$. By $\int_a^b H_\alpha\,d\alpha$ we understand the norm completion of the pre-Hilbert space that consists of continuous sections $h\in CS_i(E^\odot\odot G)$ restricted to $\RO{a,b}$ with inner product
\beqn{
\AB{h,h'}_{\SB{a,b}}
~:=~
\int_a^b\AB{h_\alpha,h'_\alpha}\,d\alpha
~=~
\int_a^b\AB{h(\alpha),h'(\alpha)}\,d\alpha.
}\eeqn
Note that all continuous sections are bounded on the compact interval $\SB{a,b}$ and, therefore, square integrable.

In order to see that (later on, in Proposition \ref{denseprop}) we will have enough sections fulfilling the stability condition in \eqref{stabsec}, it is necessary to convince ourselves that $\int_a^b H_\alpha\,d\alpha$ contains a sufficient number of sections that are only piecewise continuous. As in \cite[Proposition 4.2]{Ske07} we show:

\bprop\label{rcllprop}
$\int_a^b H_\alpha\,d\alpha$ contains the space $\eR_{\RO{a,b}}$ of restrictions to $\RO{a,b}$ of those sections $h$ for which $t\mapsto h(t)$ is right continuous with finite jumps (by this we mean, in particular, that there exists a left limit) in finitely many points of $\RO{a,b}$, and bounded on $\RO{a,b}$, as a pre-Hilbert subspace.
\eprop

(A jump at $b$ would not contribute to the inner product. So, the restriction to the right open interval $\RO{a,b}$ is necessary in order that the inner product be definite.)

Let $\sS$ denote the subspace of all sections $h=\bfam{h_t}_{t\in\R_+}$ of $E^\odot\odot G$ which are \hl{locally $\eR$}, that is, for every $0\le a<b<\infty$ the restriction of $h$ to $\RO{a,b}$ is in $\eR_{\RO{a,b}}$, and which are \hl{stable} with respect to the isometry $\zeta'_1$, that is, there exists an $\alpha_0\ge0$ such that \eqref{stabsec} holds for all $\alpha\ge\alpha_0$. By $\sN$ we denote the subspace of all sections in $\sS$ which are eventually $0$, that is, of all sections $h\in\sS$ for which there exists an $\alpha_0\ge0$ such that $h_\alpha=0$ for all $\alpha\ge\alpha_0$. A straightforward verification shows that
\beqn{
\AB{h,h'}
~:=~
\lim_{m\to\infty}\int_m^{m+1}\AB{h(\alpha),h'(\alpha)}\,d\alpha
}\eeqn
defines a semiinner product on $\sS$ and that $\AB{h,h}=0$ if and only if $h\in\sN$. Actually, we have
\beqn{
\AB{h,h'}
~=~
\int_T^{T+1}\AB{h(\alpha),h'(\alpha)}\,d\alpha
}\eeqn
for all sufficiently large $T>0$; see \cite[Lemma 2.1]{Arv06}. So, $\sS/\sN$ becomes a pre-Hilbert space with inner product $\AB{h+\sN,h'+\sN}:=\AB{h,h'}$. By $H$ we denote its completion.

As in \cite[Proposition 4.3]{Ske07} we show:

\bprop\label{denseprop}
For every section $h$ and every $\alpha_0\ge0$ define the section $h^{\alpha_0}$ as
\beqn{
h^{\alpha_0}_\alpha
~:=~
\begin{cases}
0&\alpha<\alpha_0
\\
{\zeta'_1}^nh_{\alpha-n}&\alpha\in\RO{\alpha_0+n,\alpha_0+n+1},n\in\N_0.
\end{cases}
}\eeqn
If $h$ is in $CS_i(E^\odot\odot G)$, then $h^{\alpha_0}$ is in $\sS$. Moreover, the set $\bCB{h^{\alpha_0}+\sN\colon h\in CS_i(E^\odot\odot G),\alpha_0\ge0}$ is a dense subspace of $H$.
\eprop

Observe that on $H$ we have a canonical representation of $\cB$ that acts simply pointwise on sections. (This representation is faithful, because $E^\odot$ is faithful. Nondegeneracy we see in a minute.) It is now completely plain to see that for every $t\in\R_+$ the map $x_t\odot h\mapsto x_th$, where
\beqn{
(x_th)_\alpha
~=~
\begin{cases}
x_th_{\alpha-t}&\alpha\ge t,
\\
0&\text{else},
\end{cases}
}\eeqn
defines an isometry $w_t\colon E_t\odot H\rightarrow H$, and that these isometries iterate associatively as required for a right dilation.

The following proposition may be proved as \cite[Proposition 4.6]{Ske07}. The arguments are similar to those used to show the statement about density in Proposition \ref{Hcontprop}. Actually, the proof is simpler, because thanks to Proposition \ref{rcllprop} we need not worry to obtain an approximation by continuous sections. Thus, it is not necessary to involve \it{partitions of unity}.

\bprop\label{surprop}
Each $w_t$ is surjective. In particular, for $t=0$ this means that the canonical representation of $\cB=E_0$ on $H$ is nondegenerate, so that $H$ is a correspondence from $\cB$ to $\C$.
\eprop

As in \cite[Proposition 4.6]{Ske07} we show:

\bprop\label{cdprop}
The $w_t$ are \hl{continuous} in the sense that for every continuous section $x\in CS_i(E^\odot)$ and every $h\in H$ the function $t\mapsto x_th$ is continuous.
\eprop

We are now in a position to prove the main result of these notes.

\bthm\label{rdilthm}
Let $E^\odot$ be a faithful (continuous) product system of correspondences over a \nbd{C^*}algebra $\cB$. Then $E^\odot$ admits a faithful nondegenerate (continuous) representation on a Hilbert space.
\ethm

\proof
It is clear from the preceding propositions that if $E^\odot$ is faithful and continuous, then the $w_t$ form a continuous right dilation. So, the representation $\eta^w_t$ is faithful, nondegenerate and \hl{continuous} in the sense that $t\mapsto\eta^w_t(x_t)$ is strongly continuous for every continuous section $x\in CS_i(E^\odot)$.

If $E$ is just an algebraic product system, then everything is much easier. (Simply, instead of $\int_a^b H_\alpha\,d\alpha$ use direct sums. All the technical Propositions \ref{Hcontprop}, \ref{rcllprop}, \ref{denseprop}, \ref{surprop}, \ref{cdprop} and Corollary \ref{shiftcor} are superfluous, and we obtain a right dilation that, definitely, is not continuous.)\qed

\brem\label{fufarem}
The condition to be faithful is also necessary for a product system to admit a faithful nondegenerate representation $\eta$. Indeed, suppose there is a $t\in\R_+$ such that $bx_t=0$ for all $x_t\in E_t$. Then $\eta_0(b)\eta_t(x_t)h=\eta_t(bx_t)h=0$ for all $x_t\in E_t,h\in H$. As $\eta_t$ is nondegenerate, this implies $\eta_0(b)=0$ and as $\eta$ is faithful, this implies $b=0$.

By \cite[Lemma 3.2]{Ske07} \it{continuous} product systems of correspondences over a \it{unital} \nbd{C^*}al\-ge\-bra are full automatically. Note that this is not true in the nonunital case; see \cite[Example 4.13]{Ske04p} (old version). So for having Theorem \ref{rdilthm} for correspondences over nonunital \nbd{C^*}algebras, the improved version \cite[Theorem 8.3]{Ske09} of Hirshberg's result \cite{Hir05a} is indispensable.
\erem

\section{Applications}\label{applSEC}

In this section we discuss two applications of Theorem \ref{rdilthm}. One is a \nbd{C^*}version of a theorem due to Arveson and Kishimoto, the other an existence result for a certain type of dilation of CP-semigroups.

\bthm\label{ArKiC*thm}
Let $E$ be a full Hilbert module over a unital \nbd{C^*}algebra $\cB$ and let $\vt$ be a strongly continuous strict faithful \nbd{E_0}semigroup on $\sB^a(E)$. Then there exists a faithful correspondence $K$ from $\sB^a(E)$ to $\C$ with strict left action (that is, a Hilbert space with a faithful nondegenerate strict representation of $\sB^a(E)$) and a strongly continuous unitary group $u$ on $K$ such that $\vt_t(a)k=u_tau_t^*k$ for all $a\in\sB^a(E)$, $t\in\R_+$, $k\in K$.
\ethm

\proof
Suppose we have a left dilation $v_t\colon E\odot E_t\rightarrow E$ and a right dilation $w_t\colon E_t\odot H\rightarrow H$ of a product system $E^\odot$. Then, by setting $u_t:=(v_t\odot\id_H)(\id_E\odot w_t^*)$, which acts as $x\odot y_th\mapsto xy_t\odot h$ (note that this is nontrivial!), we define a unitary semigroup $u_t$ on $K:=E\odot H$. Moreover, we recover the \nbd{E_0}semigroup $\vt^v$ on $\sB^a(E)$ by restricting $u_t\bullet u_t^*$ to $\sB^a(E)\odot\id_H$. So, our job is to recover $\vt$ as $\vt^v$ from a left dilation of a product system (mentioned in the introduction), to construct a right dilation (Theorem \ref{rdilthm}), and to show that the corresponding semigroup $u$ is sufficiently continuous.

Suppose the \nbd{E_0}semigroup $\vt$ acts on a Hilbert \nbd{\cB}module $E$ with a unit vector $\xi$. (Otherwise, by \cite[Lemma 3.2]{Ske09} there is a natural number $n$ such that $E^n$ has a unit vector. We may pass to the inflation of $\vt$ to an \nbd{E_0}semigroup on $\sB^a(E^n)=M_n(\sB^a(E))$ that gives back $\vt$ by embedding $\sB^a(E)$  unitally into the diagonal of $M_n(\sB^a(E))$.) Then as explained in \cite[Section 7]{Ske03b}, we obtain a continuous product system $E^\odot$ and a left dilation of $E^\odot$ that gives back $\vt$, in the following way. Put $E_t:=\vt_t(\xi\xi^*)E$, equip $E_t$ with the left action $bx_t:=\vt_t(\xi b\xi^*)x_t$ of $\cB$, and define $v_t\colon x\odot y_t\mapsto\vt_t(x\xi^*)y_t$. Then $x_sy_t:=v_t(x_s\odot y_t)$ defines a product system structure on $E^\odot=\bfam{E_t}_{t\in\R_+}$ and the $v_t$ define a left dilation of $E^\odot$ to $E$ giving back $\vt=\vt^v$. Moreover, choosing for $i_t$ the canonical embedding of the submodule $E_t$ of $E$ into $E$, we turn $E^\odot$ into a continuous product system. Note that the left dilation is \hl{continuous} in the sense that $t\mapsto xy_t=\vt_t(x\xi^*)y_t=\vt_t(x\xi^*)i_ty_t$ for every continuous section $y$.

Since $E$ is full and $\vt_t$ is faithful, $E_t$ is faithful. For $w_t$ we choose the right dilation from Theorem \ref{rdilthm}.

The proof of continuity of $u_t$ is quite standard following similar proofs in \cite{Ske03b,Ske07}: As the family $u_t$ is bounded uniformly, it is sufficient to check strong continuity on the total subset of $K$ formed by all $x\odot h$. Let $\zeta$ denote a continuous section of unit vectors with $\zeta_0=\U$, as granted by \cite[Lemma 3.2]{Ske07}. Then for all sufficiently small $t$, the elements $\zeta_th$ are close to $h$ and, therefore, the elements $x\odot\zeta_th$ are close to $x\odot h$. Similarly, $x\zeta_t$ is close to $x$ and, therefore, also $x\zeta_t\odot h$ is close to $x\odot h$. We find that
\beqn{
u_t(x\odot h)-x\odot h
~=~
u_t(x\odot h)-u_t(x\odot\zeta_th)+u_t(x\odot\zeta_th)-x\odot h
~=~
u_t(x\odot h-x\odot\zeta_th)+(x\zeta_t\odot h-x\odot h)
}\eeqn
is close to $0$.\qed

\bcor\label{ArKiC*cor}
Let $\vt$ be a strongly continuous faithful semigroup of nondegenerate endomorphisms on a \nbd{C^*}algebra $\cB$. (By nondegenerate we mean $\cls\vt_t(\cB)\cB=\cB$ for all $t\in\R_+$.) Then there exists a faithful correspondence $K$ from $\cB$ to $\C$ and a strongly continuous unitary group $u$ on $K$ such that $\vt_t(b)k=u_tbu_t^*k$ for all $b\in\cB$, $t\in\R_+$, $k\in K$.
\ecor

\proof
$\vt$ extends from $\cB=\sK(\cB)$ to a strict \nbd{E_0}semigroup on $\sB^a(\cB)$.\qed

\brem
For von Neumann algebras $\cB\subset\sB(G)$ and \nbd{E_0}semigroups that are continuous in the strong operator topology of $\sB(G)$, the statement of Corollary \ref{ArKiC*cor} is due to Arveson and Kishimoto \cite{ArKi92}. In \cite{Ske07p} we will provide analogue proofs of a version of Theorem \ref{ArKiC*cor} for von Neumann modules and the result of \cite{ArKi92}, that is, Corollary \ref{ArKiC*cor} for von Neumann algebras.
\erem

We now come to \it{elementary} dilations of CP-semigroups. Recall that a \hl{CP-semigroup} is a semigroup $T=\bfam{T_t}_{t\in\R_+}$ of completely positive (CP-)maps on a \nbd{C^*}algebra. A CP-map $T$ is \hl{faithful}, if $T(b^*b)=0$ implies $b=0$ for all $b\in\cB$. It is \hl{semifaithful}, if its GNS-correspondence $E$ is faithful. (Recall that by Paschke \cite{Pas73} the GNS-correspondence is that unique correspondence over $\cB$ that contains a vector $\xi$ that generates $E$ as correspondence and that fulfills $\AB{\xi,b\xi}=T(b)$. Alternatively one may require that some (and, therefore, every) Stinespring representation is faithful.) A CP-semigroup $T$ is \hl{(semi-)faithful}, if every $T_t$ is (semi-)faithful. A CP-semigroup $T$ is \hl{elementary}, if it has the form $T_t(b)=c_t^*bc_t$ for some semigroup $c=\bfam{c_t}_{t\in\R_+}$ of elements in $\cB$. An \hl{elementary dilation} of a CP-semigroup $T$ on $\cB$ is a \nbd{C^*}algebra $\cA$ with an embedding $\vp\colon\cB\rightarrow\cA$ and a semigroup $c=\bfam{c_t}_{t\in\R_+}$ of elements in $\cA$ such that
\beqn{
\vp\circ T(b)
~=~
c^*_t\vp(b)c_t
}\eeqn
for all $b\in\cB$ and $t\in\R_+$.

\bthm\label{elemdilthm}
Every semifaithful (strongly continuous) CP-semigroup on a unital \nbd{C^*}algebra admits a (strongly continuous) elementary dilation to some $\sB(H)$.
\ethm

\proof
Recall that a \hl{unit} for a product system $E^\odot$ of correspondences $E_t$ over a unital \nbd{C^*}al\-ge\-bra $\cB$ is a family $\xi^\odot=\bfam{\xi_t}$ of elements $\xi_t\in E_t$ that factors as $\xi_{s+t}=\xi_s\xi_t$ with $\xi_0=\U$. Bhat and Skeide \cite{BhSk00} associate with every CP-semigroup $T$ on a unital \nbd{C^*}algebra $\cB$ a product system $E^\odot$ of correspondences over $\cB$ and a unit $\xi^\odot$ such that $T_t(b)=\AB{\xi_t,b\xi_t}$. (This product system is unique, if we require that it is generated by the unit $\xi^\odot$.)

If $T$ is semifaithful, then already the left action of $\cB$ on the \nbd{\cB}bimodule of $E_t$ generated by $\xi_t$ is faithful. Therefore, we may apply Theorem \ref{rdilthm} to obtain a right dilation $w$ of $E^\odot$ to a correspondence $H$ from $\cB$ to $\C$ and, further, e faithful nondegenerate representation $\eta:=\eta^w$ of $E^\odot$ on $H$. It follows that $c_t:=\eta_t(\xi_t)$ is a semigroup in $\sB(H)$ and that $c_t^*\eta_0(b)c_t=\eta_0(\AB{\xi_t,b\xi_t})=\eta_0\circ T_t(b)$.

If $T$ is strongly continuous then by \cite[Section 7]{Ske03b} the product system generated by $\xi^\odot$ is continuous and $\xi^\odot$ is among the continuous sections. By Theorem \ref{rdilthm}, we may chose the right dilation $w$ continuous. But this means precisely that $t\mapsto\eta_t(x_t)$ is strongly continuous for every continuous section $x$. Since $\xi^\odot$ is continuous, so is $c$.\qed

\brem
Note that, in particular, uniformly continuous CP-semigroups fit into the assumptions of Theorem \ref{elemdilthm}. (In fact, every $T_t$ is invertible so that $b\xi_t=0\Rightarrow T_t(b^*b)=0\Rightarrow b^*b=0\Rightarrow b=0$, so that $T_t$ is even faithful.) We think that in this way Theorem \ref{elemdilthm} might be helpful in finding a proof of the fact that a suitable strong closure of $E^\odot$ contains a \hl{central} continuous unit $\om^\odot$, that is, a unit where all $\om_t$ commute with all elements of $\cB$. By Barreto, Bhat, Liebscher and Skeide \cite{BBLS04} this statement is equivalent to the results by Christensen and Evans \cite{ChrEv79} on the form of the generator of a uniformly continuous CP-semigroup.
\erem

\brem
For a CP-semigroup $T$ on $\cB$, the proof of Theorem \ref{elemdilthm} starts with the construction from \cite{BhSk00} of a product system $E^\odot$ of \nbd{\cB}correspondences and a unit $\xi^\odot$ for that product system. Then Theorem \ref{rdilthm} is applied to that product system $E^\odot$ providing a representation of $E^\odot$. In the case of a normal CP-semigroup on a von Neumann algebra $\cB$, also Muhly and Solel \cite{MuSo02} construct a product system and, then, a representation of that product system. We would like to emphasize that the product system constructed in \cite{MuSo02} is $E'^\odot$, the commutant of $E^\odot$ (see \cite{Ske03c,Ske08}). The construction in \cite{MuSo02} of a representation (that is, of a right dilation) of $E'^\odot$ is equivalent (via commutant) to the older construction in \cite{BhSk00} of an \nbd{E_0}semigroup (that is, a left dilation) for $E^\odot$ based on the unit $\xi^\odot$. The construction of an \nbd{E_0}semigroup from a unit is much easier than the general result in \cite{Ske07} without unit. It has nothing to do with our representation in Theorem \ref{rdilthm}.
\erem

\brem
We would like to mention that, after the first version of these notes, in Skeide \cite{Ske09a} we constructed a continuous product system with a continuous left dilation for a strongly continuous \nbd{E_0}semigroup on $\sB^a(E)$, where $E$ may be a full Hilbert module over a \nbd{\sigma}unital \nbd{C^*}algebra $\cB$, or any Hilbert module $E$ over an arbitrary \nbd{C^*}algebra as long as $E^\infty$ contains a direct summand $\cB$. Theorem \ref{ArKiC*thm} holds under these conditions, too.
\erem

\setlength{\baselineskip}{2.5ex}

\newcommand{\Swap}[2]{#2#1}\newcommand{\Sort}[1]{}
\providecommand{\bysame}{\leavevmode\hbox to3em{\hrulefill}\thinspace}
\providecommand{\MR}{\relax\ifhmode\unskip\space\fi MR }
\providecommand{\MRhref}[2]{%
  \href{http://www.ams.org/mathscinet-getitem?mr=#1}{#2}
}
\providecommand{\href}[2]{#2}


\end{document}